\documentclass{article}
\usepackage{spconf,amsmath,amssymb,amsthm,graphicx,stmaryrd}

\newtheorem{theorem}{Theorem}
\newtheorem{proposition}{Proposition}

\title{SINGULAR VALUES AND EIGENVALUES OF TENSORS: A VARIATIONAL APPROACH}
\name{Lek-Heng~Lim\thanks{This work appeared in: \textit{Proceedings
of the IEEE International Workshop on Computational Advances in
Multi-Sensor Adaptive Processing} (CAMSAP '05), \textbf{1} (2005),
pp.\ 129--132.}}
\address{Stanford University\\
Institute for Computational and Mathematical Engineering\\
Gates Building 2B, Room 286, Stanford, CA 94305}

\begin{document}
\ninept

\maketitle

\begin{abstract}
We propose a theory of eigenvalues, eigenvectors, singular values,
and singular vectors for tensors based on a constrained variational
approach much like the Rayleigh quotient for symmetric matrix
eigenvalues. These notions are particularly useful in generalizing
certain areas where the spectral theory of matrices has
traditionally played an important role. For illustration, we will
discuss a multilinear generalization of the Perron-Frobenius
theorem.
\end{abstract}

\section{INTRODUCTION}

It is well known that the eigenvalues and eigenvectors of a symmetric matrix
$A$ are the critical values and critical points of its Rayleigh quotient,
$\mathbf{x}^{\intercal}A\mathbf{x}/\lVert\mathbf{x}\rVert_{2}^{2}$, or
equivalently, the critical values and points of the quadratic form
$\mathbf{x}^{\intercal}A\mathbf{x}$ constrained to vectors with unit $l^{2}$-norm,
$\{\mathbf{x}\mid\lVert\mathbf{x}\rVert_{2}=1\}$. If $L:\mathbb{R}^{n}%
\times\mathbb{R\rightarrow R}$ is the associated Lagrangian with Lagrange
multiplier $\lambda$,%
\[
\quad L(\mathbf{x},\lambda)=\mathbf{x}^{\intercal}A\mathbf{x}-\lambda
(\lVert\mathbf{x}\rVert_{2}^{2}-1),
\]
then the vanishing of $\nabla L$ at a critical point $(\mathbf{x}_{c}%
,\lambda_{c})\in\mathbb{R}^{n}\times\mathbb{R}$ yields the familiar defining
condition for eigenpairs%
\begin{equation}
A\mathbf{x}_{c}=\lambda_{c}\mathbf{x}_{c}. \label{sevd}%
\end{equation}
Note that this approach does not work if $A$ is nonsymmetric --- the critical
points of $L$ would in general be different from the solutions of
\eqref{sevd}%
.

A little less widely known is an analogous variational approach to the
singular values and singular vectors of a matrix $A\in\mathbb{R}^{m\times n}$,
with $\mathbf{x}^{\intercal}A\mathbf{y}/\mathbf{\lVert\mathbf{x}\rVert}%
_{2}\mathbf{\lVert\mathbf{y}\rVert}_{2}$ assuming the role of the Rayleigh
quotient. The associated Lagrangian function $L:\mathbb{R}^{m}\times
\mathbb{R}^{n}\times\mathbb{R}\rightarrow\mathbb{R}$ is now%
\[
L(\mathbf{x},\mathbf{y},\sigma)=\mathbf{x}^{\intercal}A\mathbf{y}-\sigma
(\lVert\mathbf{x}\rVert_{2}\lVert\mathbf{y}\rVert_{2}-1).
\]
$L$ is continuously differentiable for non-zero $\mathbf{x},\mathbf{y}$. The
first order condition yields%
\[
A\mathbf{y}_{c}/\lVert\mathbf{y}_{c}\rVert_{2}=\sigma_{c}\mathbf{x}_{c}%
/\lVert\mathbf{x}_{c}\rVert_{2},\quad A^{\intercal}\mathbf{x}_{c}/\lVert
\mathbf{x}_{c}\rVert_{2}=\sigma_{c}\mathbf{y}_{c}/\lVert\mathbf{y}_{c}%
\rVert_{2},
\]
at a critical point $(\mathbf{x}_{c},\mathbf{y}_{c},\sigma_{c})\in
\mathbb{R}^{m}\times\mathbb{R}^{n}\times\mathbb{R}$. Writing $\mathbf{u}%
_{c}=\mathbf{x}_{c}/\lVert\mathbf{x}_{c}\rVert_{2}$ and $\mathbf{v}%
_{c}=\mathbf{y}_{c}/\lVert\mathbf{y}_{c}\rVert_{2}$, we get the familiar%
\begin{equation}
A\mathbf{v}_{c}=\sigma_{c}\mathbf{u}_{c},\quad A^{\intercal}\mathbf{u}_{c}%
=\sigma_{c}\mathbf{v}_{c}. \label{svd}%
\end{equation}

Although it is not immediately clear how the usual definitions of eigenvalues
and singular values via
\eqref{sevd}
and
\eqref{svd}
may be generalized to tensors of order $k\geq3$ (a matrix is regarded as an
order-$2$ tensor), the constrained variational approach generalizes in a
straight-forward manner --- one simply replaces the bilinear functional
$\mathbf{x}^{\intercal}A\mathbf{y}$ (resp.\ quadratic form $\mathbf{x}^{\intercal
}A\mathbf{x}$) by the multilinear functional (resp.\ homogeneous polynomial)
associated with a tensor (resp.\ symmetric tensor) of order $k$. The
constrained critical values/points then yield a notion of singular
values/vectors (resp.\ eigenvalues/vectors) for order-$k$ tensors.

An important point of distinction between the order-$2$ and order-$k$ cases is
in the choice of norm for the constraints. At first glance, it may appear that
we should retain the $l^{2}$-norm. However, the criticality conditions so
obtained are no longer \textit{scale invariant} (ie.\ the property that
$\mathbf{x}_{c}$ in
\eqref{sevd}
or $(\mathbf{u}_{c},\mathbf{v}_{c})$ in
\eqref{svd}
may be replaced by $\alpha\mathbf{x}_{c}$ or $(\alpha\mathbf{u}_{c}%
,\alpha\mathbf{v}_{c})$ without affecting the validity of the equations). To
preserve the scale invariance of eigenvectors and singular vectors for tensors
of order $k\geq3$, the $l^{2}$-norm must be replaced by the $l^{k}$-norm
(where $k$ is the order of the tensor),%
\[
\lVert\mathbf{x}\rVert_{k}=(\lvert x_{1}\rvert^{k}+\dots+\lvert x_{n}%
\rvert^{k})^{1/k}.
\]
The consideration of eigenvalues and singular values with respect to $l^{p}%
$-norms where $p\neq2$ is prompted by recent works \cite{Cho1, Cho2} of
Choulakian, who studied such notions for matrices.

Nevertheless, we shall not insist on having scale invariance. Instead, we will
define eigenpairs and singular pairs of tensors with respect to any $l^{p}%
$-norm ($p>1$) as they can be interesting even when $p\neq k$. For example,
when $p=2$, our defining equations for singular values/vectors
\eqref{msvd}
become the equations obtained in the best rank-$1$ approximations of tensors
studied by Comon \cite{Comon} and de Lathauwer et.\ al.\ \cite{dLMV2}. For the
special case of symmetric tensors, our equations for eigenvalues/vectors for
$p=2$ and $p=k$ define respectively, the Z-eigenvalues/vectors and
H-eigenvalues/vectors in the soon-to-appear paper \cite{Q} of Qi. For
simplicity, we will restrict our study to integer-valued $p$ in this paper.

We thank Gunnar Carlsson, Pierre Comon, Lieven de Lathauwer, Vin de Silva, and
Gene Golub for helpful discussions. We would also like to thank Liqun Qi for
sending us an advanced copy of his very relevant preprint.

\section{TENSORS AND MULTILINEAR FUNCTIONALS}

A $k$-array of real numbers representing an order-$k$ tensor will be denoted
by $A=%
\llbracket
a_{j_{1}\cdots j_{k}}%
\rrbracket
\in\mathbb{R}^{d_{1}\times\dots\times d_{k}}$. Just as an order-$2$ tensor
(ie.\ matrix) may be multiplied on the left and right by a pair of matrices
(of consistent dimensions), an order-$k$ tensor may be `multiplied on $k$
sides' by $k$ matrices. The \textit{covariant multilinear matrix
multiplication} of $A$ by matrices $M_{1}=[m_{j_{1}i_{1}}^{(1)}]\in
\mathbb{R}^{d_{1}\times s_{1}},\dots,M_{k}=[m_{j_{k}i_{k}}^{(k)}]\in
\mathbb{R}^{d_{k}\times s_{k}}$ is defined by%
\begin{multline*}
A(M_{1},\dots,M_{k}):=\\%
\Bigl
\llbracket
\sum\nolimits_{j_{1}=1}^{d_{1}}\cdots\sum\nolimits_{j_{k}=1}^{d_{k}}%
a_{j_{1}\cdots j_{k}}m_{j_{1}i_{1}}^{(1)}\cdots m_{j_{k}i_{k}}^{(k)}%
\Bigr
\rrbracket
\in\mathbb{R}^{s_{1}\times\dots\times s_{k}}.
\end{multline*}
This operation arises from the way a multilinear functional transforms under
compositions with linear maps. In particular, the multilinear functional
associated with a tensor $A\in\mathbb{R}^{d_{1}\times\dots\times d_{k}}$ and
its gradient may be succinctly expressed via covariant multilinear
multiplication:%
\begin{gather}
A(\mathbf{x}_{1},\dots,\mathbf{x}_{k})=\sum\nolimits_{j_{1}=1}^{d_{1}}%
\cdots\sum\nolimits_{j_{k}=1}^{d_{k}}a_{j_{1}\cdots j_{k}}x_{j_{1}}%
^{(1)}\cdots x_{j_{k}}^{(k)},\label{mmm}\\
\nabla_{\mathbf{x}_{i}}A(\mathbf{x}_{1},\dots,\mathbf{x}_{k})=A(\mathbf{x}%
_{1},\dots,\mathbf{x}_{i-1},I_{d_{i}},\mathbf{x}_{i+1},\dots,\mathbf{x}%
_{k}).\nonumber
\end{gather}
Note that we have slightly abused notations by using $A$ to denote both the
tensor and its associated multilinear functional.

An order-$k$ tensor $%
\llbracket
a_{j_{1}\cdots j_{k}}%
\rrbracket
\in\mathbb{R}^{n\times\dots\times n}$ is called \textit{symmetric}
if $a_{j_{\sigma(1)}\cdots j_{\sigma(k)}}=a_{j_{1}\cdots j_{k}}$ for
any permutation $\sigma\in\mathfrak{S}_{k}$. The homogeneous
polynomial associated with
a symmetric tensor $A=%
\llbracket
a_{j_{1}\cdots j_{k}}%
\rrbracket
$ and its gradient can again be conveniently expressed as%
\begin{align}
A(\mathbf{x},\dots,\mathbf{x})  &  =\sum\nolimits_{j_{1}=1}^{n}\cdots
\sum\nolimits_{j_{k}=1}^{n}a_{j_{1}\cdots j_{k}}x_{j_{1}}\cdots x_{j_{k}%
},\label{sm}\\
\nabla A(\mathbf{x},\dots,\mathbf{x})  &  =kA(I_{n},\mathbf{x},\dots
,\mathbf{x}).\nonumber
\end{align}
Observe that for a symmetric tensor $A$,
\begin{multline}
A(I_{n},\mathbf{x},\mathbf{x},\dots,\mathbf{x})=A(\mathbf{x},I_{n}%
,\mathbf{x},\dots,\mathbf{x})=\\
\dots=A(\mathbf{x},\mathbf{x},\dots,\mathbf{x},I_{n}). \label{vecdir}%
\end{multline}

The preceding discussion is entirely algebraic but we will now introduce
norms on the respective spaces. Let $\lVert\cdot\rVert_{\alpha_{i}}$ be a
norm on $\mathbb{R}^{d_{i}}$, $i=1,\dots,k$. Then the \textit{norm} (cf.\
\cite{DF}) of the multilinear functional
$A:\mathbb{R}^{d_{1}}\times\dots\times
\mathbb{R}^{d_{k}}\rightarrow\mathbb{R}$ induced by $\lVert\cdot\rVert
_{\alpha_{1}},\dots,\lVert\cdot\rVert_{\alpha_{k}}$ is defined as%
\[
\lVert A\rVert_{\alpha_{1},\dots,\alpha_{k}}:=\sup\frac{\lvert A(\mathbf{x}%
_{1},\dots,\mathbf{x}_{k})\rvert}{\lVert\mathbf{x}_{1}\rVert_{\alpha_{1}%
}\cdots\lVert\mathbf{x}_{k}\rVert_{\alpha_{k}}}%
\]
where the supremum is taken over all non-zero $\mathbf{x}_{i}\in
\mathbb{R}^{d_{i}}$, $i=1,\dots,k$. We will be interested in the
case where the $\lVert\cdot\rVert_{\alpha_{i}}$'s are $l^{p}$-norms.
Recall that for $1 \le p \le \infty$, the $l^{p}$-norm is a
continuously differentiable function on
$\mathbb{R}^{n}\backslash\{\mathbf{0}\}$. For
$\mathbf{x}=[x_{1},\dots,x_{n}]^{\intercal}\in\mathbb{R}^{n}$, we will write%
\[
\mathbf{x}^{p}:=[x_{1}^{p},\dots,x_{n}^{p}]^{\intercal}%
\]
(ie.\ taking $p$th power coordinatewise) and%
\[
\varphi_{p}(\mathbf{x):=}[\operatorname{sgn}(x_{1})x_{1}^{p},\dots
,\operatorname{sgn}(x_{n})x_{n}^{p}]^{\intercal}%
\]
where%
\[
\operatorname{sgn}(x)=%
\begin{cases}
+1 & \text{if }x>0,\\
0 & \text{if }x=0,\\
-1 & \text{if }x<0.
\end{cases}
\]
Observe that if $p$ is even, then $\varphi_{p}(\mathbf{x})=\mathbf{x}^{p}$.
The gradient of the $l^{p}$-norm is given by%
\[
\nabla\lVert\mathbf{x}\rVert_{p}=\frac{\varphi_{p-1}(\mathbf{x)}}%
{\lVert\mathbf{x}\rVert_{p}^{p-1}}%
\]
or simply $\nabla\lVert\mathbf{x}\rVert_{p}=\mathbf{x}^{p-1}/\lVert
\mathbf{x}\rVert_{p}^{p-1}$ when $p$ is even.

\section{SINGULAR VALUES AND SINGULAR VECTORS}

Let $A\in\mathbb{R}^{d_{1}\times\dots\times d_{k}}$. Then $A$ defines a
multilinear functional $A:\mathbb{R}^{d_{1}}\times\dots\times\mathbb{R}%
^{d_{k}}\rightarrow\mathbb{R}$ via
\eqref{mmm}%
. Let us equip $\mathbb{R}^{d_{i}}$ with the $l^{p_{i}}$-norm, $\lVert
\cdot\rVert_{p_{i}}$, $i=1,\dots,k$. We will define the singular values and
singular vectors of $A$ as the critical values and critical points of
$A(\mathbf{x}_{1},\dots,\mathbf{x}_{k})/\lVert\mathbf{x}_{1}\rVert_{p_{1}%
}\cdots\lVert\mathbf{x}_{k}\rVert_{p_{k}}$, suitably normalized. Taking a
constrained variational approach, we let $L:\mathbb{R}^{d_{1}}\times
\dots\times\mathbb{R}^{d_{k}}\times\mathbb{R}\rightarrow\mathbb{R}$ be%
\begin{multline*}
L(\mathbf{x}_{1},\dots,\mathbf{x}_{k},\sigma) :=\\
A(\mathbf{x}_{1},\dots,\mathbf{x}_{k})-\sigma(\lVert\mathbf{x}_{1}%
\rVert_{p_{1}}\cdots\lVert\mathbf{x}_{k}\rVert_{p_{k}}-1).
\end{multline*}
$L$ is continuously differentiable when $\mathbf{x}_{i}\neq\mathbf{0}$,
$i=1,\dots,k$. The vanishing of the gradient,%
\[
\nabla L=(\nabla_{\mathbf{x}_{1}}L,\dots,\nabla_{\mathbf{x}_{k}}%
L,\nabla_{\sigma}L)=(\mathbf{0},\dots,\mathbf{0},0)
\]
gives
\begin{equation}\label{msvd}
\begin{aligned}
A(I_{d_{1}},\mathbf{x}_{2},\mathbf{x}_{3},\dots,\mathbf{x}_{k}) &
=\sigma\varphi_{p_{1}-1}(\mathbf{x}_{1}),\\
A(\mathbf{x}_{1},I_{d_{2}},\mathbf{x}_{3},\dots,\mathbf{x}_{k}) &
=\sigma\varphi_{p_{2}-1}(\mathbf{x}_{2}),\\
& \vdots\\
A(\mathbf{x}_{1},\mathbf{x}_{2},\dots,\mathbf{x}_{k-1},I_{d_{k}}) &
=\sigma\varphi_{p_{k}-1}(\mathbf{x}_{k}), \end{aligned}
\end{equation}
at a critical point $(\mathbf{x}_{1},\dots,\mathbf{x}_{k},\sigma)$. As in the
derivation of
\eqref{svd}%
, one gets also the unit norm condition
\[
\lVert\mathbf{x}_{1}\rVert_{p_{1}}=\cdots=\lVert\mathbf{x}_{k}\rVert_{p_{k}%
}=1.
\]
The unit vector $\mathbf{x}_{i}$ and $\sigma$ in
\eqref{msvd}%
, will be called the \textit{mode-}$i$\textit{\ singular vector},
$i=1,\dots,k$, and \textit{singular value} of $A$ respectively. Note that the
mode-$i$ singular vectors are simply the order-$k$ equivalent of left- and
right-singular vectors for order $2$ (a matrix has two `sides' or modes while
an order-$k$ tensor has $k$).

We will use the name $l^{p_{1},\dots,p_{k}}$-singular values/vectors if we
wish to emphasize the dependence of these notions on $\lVert\cdot\rVert
_{p_{i}}$, $i=1,\dots,k$. If $p_{1}=\dots=p_{k}=p$, then we will use the
shorter name $l^{p}$-singular values/vectors. Two particular choices of $p$
will be of interest to us: $p=2$ and $p=k$ --- both of which reduce to the
matrix case when $k=2$ (not so for other choices of $p$). The former yields
\[
A(\mathbf{x}_{1},\dots,\mathbf{x}_{i},I_{d_{i}},\mathbf{x}_{i+1}%
,\dots,\mathbf{x}_{k})=\sigma\mathbf{x}_{i},\quad i=1,\dots,k,
\]
while the latter yields a homogeneous system of equations that is invariant
under scaling of $(\mathbf{x}_{1},\dots,\mathbf{x}_{k})$. In fact, when $k$ is
even, the $l^{p}$-singular values/vectors are solutions to%
\[
A(\mathbf{x}_{1},\dots,\mathbf{x}_{i},I_{d_{i}},\mathbf{x}_{i+1}%
,\dots,\mathbf{x}_{k})=\sigma\mathbf{x}_{i}^{k-1},\quad i=1,\dots,k.
\]

The following results are easy to show. The first proposition follows from the
definition of norm and the observation that a maximizer in an open set must be
critical. The second proposition follows from the definition of
hyperdeterminant \cite{GKZ2}; the conditions on $d_{i}$ are necessary and
sufficient for the existence of $\Delta$.

\begin{proposition}
The largest $l^{p_{1},\dots,p_{k}}$-singular value is equal to the norm of the
multilinear functional associated with $A$ induced by the norms $\lVert
\cdot\rVert_{p_{1}},\dots,\lVert\cdot\rVert_{p_{k}}$, ie.
\[
\sigma_{\max}(A) = \lVert A\rVert_{p_{1},\dots,p_{k}}.
\]

\end{proposition}

\begin{proposition}
Let $d_{1},\dots,d_{k}$ be such that%
\[
d_{i}-1\leq\sum\nolimits_{j\neq i}(d_{j}-1)\quad\text{for all }i=1,\dots,k,
\]
and $\Delta$ denote the hyperdeterminant in $\mathbb{R}^{d_{1}\times
\dots\times d_{k}}$. Then $0$ is an $l^{2}$-singular value of $A\in
\mathbb{R}^{d_{1}\times\dots\times d_{k}}$ if and only if%
\[
\Delta(A)=0.
\]

\end{proposition}

\section{EIGENVALUES AND EIGENVECTORS OF SYMMETRIC TENSORS}

Let $A\in\mathbb{R}^{n\times\dots\times n}$ be an order-$k$ \textit{symmetric}
tensor. Then $A$ defines a degree-$k$ homogeneous polynomial function
$A:\mathbb{R}^{n}\rightarrow\mathbb{R}$ via
\eqref{sm}%
. With a choice of $l^{p}$-norm on $\mathbb{R}^{n}$, we may consider the
multilinear Rayleigh quotient $A(\mathbf{x},\dots,\mathbf{x})/\lVert
\mathbf{x}\rVert_{p}^{k}$. The Lagrangian $L:\mathbb{R}^{n}\times
\mathbb{R\rightarrow R}$,%
\[
L(\mathbf{x},\lambda):=A(\mathbf{x},\dots,\mathbf{x})-\lambda(\lVert
\mathbf{x}\rVert_{p}^{k}-1),
\]
is continuously differentiable when $\mathbf{x}\neq\mathbf{0}$ and%
\[
\nabla L=(\nabla_{\mathbf{x}}L,\nabla_{\lambda}L)=(\mathbf{0},0)
\]
gives%
\begin{equation}
A(I_{n},\mathbf{x},\dots,\mathbf{x})=\lambda\varphi_{p-1}(\mathbf{x})
\label{msevd}%
\end{equation}
at a critical point $(\mathbf{x},\lambda)$ where $\lVert\mathbf{x}\rVert
_{p}=1$. The unit vector $\mathbf{x}$ and scalar $\lambda$ will be called an
$l^{p}$\textit{-eigenvector} and $l^{p}$\textit{-eigenvalue} of $A$
respectively. Note that the \textsc{lhs} in
\eqref{msevd}
satisfies the symmetry in
\eqref{vecdir}%
.

As in the case of singular values/vectors, the instances where $p=2$ and $p=k
$ are of particular interest. The $l^{2}$-eigenpairs are characterized by%
\[
A(I_{n},\mathbf{x},\dots,\mathbf{x})=\lambda\mathbf{x}%
\]
where $\lVert\mathbf{x}\rVert_{2}=1$. When the order $k$ is even, the $l^{k}%
$-eigenpairs are characterized by%
\begin{equation}
A(I_{n},\mathbf{x},\dots,\mathbf{x})=\lambda\mathbf{x}^{k-1} \label{kmsevd}%
\end{equation}
and in this case the unit-norm constraint is superfluous since
\eqref{kmsevd}
is a homogeneous system and $\mathbf{x}$ may be scaled by any non-zero scalar
$\alpha$.

We shall refer the reader to \cite{Q} for some interesting results on $l^{2}$-eigenvalues
and $l^{k}$-eigenvalues for symmetric tensors --- many of which mirrors
familiar properties of matrix eigenvalues.

\section{EIGENVALUES AND EIGENVECTORS OF NONSYMMETRIC TENSORS}

We know that one cannot use the variational approach to characterize
eigenvalues/vectors of \textit{nonsymmetric} matrices. So for an
\textit{nonsymmetric} tensor $A\in\mathbb{R}^{n\times\dots\times n}$, we will
instead \textit{define} eigenvalues/vectors by
\eqref{msevd}%
\ --- an approach that is consistent with the matrix case. As
\eqref{vecdir}
no longer holds, we now have $k$ different forms of
\eqref{msevd}%
:%
\begin{equation}
\label{mevd}\begin{aligned}
A(I_{n},\mathbf{x}_{1},\mathbf{x}_{1},\dots,\mathbf{x}_{1}) & =\mu
_{1}\varphi_{p-1}(\mathbf{x}_{1}),\\
A(\mathbf{x}_{1},I_{n},\mathbf{x}_{2},\dots,\mathbf{x}_{2}) & =\mu
_{2}\varphi_{p-1}(\mathbf{x}_{2}),\\ & \vdots\\
A(\mathbf{x}_{k},\mathbf{x}_{k},\dots,\mathbf{x}_{k},I_{n}) & =\mu
_{k}\varphi_{p-1}(\mathbf{x}_{k}). \end{aligned}
\end{equation}
We will call the unit vector $\mathbf{x}_{i}$ a mode-$i$
\textit{eigenvector} of $A$ corresponding to the mode-$i$
\textit{eigenvalue} $\mu_{i}$, $i=1,\dots,k$. Note that these are nothing
more than the order-$k$ equivalent of left and right eigenvectors.

\section{APPLICATIONS}

Several distinct generalizations of singular values/vectors and
eigenvalues/vectors from matrices to higher-order tensors have been
proposed in \cite{Comon, dLMV1, dLMV2, RC, Q}. As one can expect,
there is no one single generalization that preserves all properties
of matrix singular values/vectors or matrix eigenvalues/vectors.
In the lack of a canonical generalization, the validity of
a multilinear generalization of a bilinear
concept is often measured by the extent to which it may be applied to obtain interesting or
useful results.

The proposed notions of $l^2$- and $l^k$-singular/eigenvalues arise
naturally in the context of several different applications. We have
mentioned the relation between the $l^{2}$-singular values/vectors and the
best rank-$1$ approximation of a tensor under the Frobenius norm obtained in
\cite{Comon, dLMV2}. Another example is the appearance of
$l^{2}$-eigenvalues/vectors of symmetric tensors in the approximate
solutions of constraint satisfaction problems \cite{ADKK, DKKV}. A third
example is the use of $l^{k}$-eigenvalues for order-$k$ symmetric tensors
($k$ even) for characterizing the positive definiteness of homogeneous
polynomial forms --- a problem that is important in automatic control and
array signal processing (see \cite{Q, QT} and the references cited therein).

Here we will give an application of $l^{k}$-eigenvalues and eigenvectors of a
\textit{nonsymmetric} tensor of order $k$. We will show that a multilinear
generalization of the Perron-Frobenius theorem \cite{BP} may be deduced from
the notion of $l^{k}$-eigenvalues/vectors as defined by
\eqref{mevd}%
.

Let $A=%
\llbracket
a_{j_{1}\cdots j_{k}}%
\rrbracket
\in\mathbb{R}^{n\times\dots\times n}$. We write $A>0$ if all $a_{j_{1}\cdots
j_{k}}>0$ (likewise for $A\geq0$). We write $A>B$ if $A-B>0$ (likewise for
$A\geq B$).

An order-$k$ tensor $A$ is \textit{reducible} if there exists a
permutation $\sigma\in\mathfrak{S}_{n}$ such that the permuted
tensor
\[%
\llbracket
b_{i_{1}\cdots i_{k}}%
\rrbracket
=%
\llbracket
a_{\sigma(j_{1})\cdots\sigma(j_{k})}%
\rrbracket
\in\mathbb{R}^{n\times\dots\times n}
\]
has the property that for some $m \in\{1,\dots, n-1\}$, $b_{i_{1}\cdots i_{k}%
}=0$ for all $i_{1}\in\{1,\dots,n-m\}$ and all $i_{2},\dots,i_{k}\in
\{1,\dots,m\}$.

If we allow a few analogous matrix terminologies, then $A$ is
reducible if there exists a permutation matrix $P$ so that
\[
B=A(P,\dots,P)\in\mathbb{R}^{n\times\dots\times n}
\]
can be partitioned into $2^{n}$ subblocks and regarded as a $2\times
\dots\times2$ block-tensor with `square diagonal blocks' $B_{00\cdots0}%
\in\mathbb{R}^{m\times\dots\times m}$, $B_{11\cdots1}\in\mathbb{R}%
^{(n-m)\times\dots\times(n-m)}$, and a zero `corner block' $B_{10\cdots0}%
\in\mathbb{R}^{(n-m)\times m\times\dots\times m}$ which we may assume without
loss of generality to be in the $(1,0,\dots,0)$-`corner'.

We say that $A$ is \textit{irreducible} if it is not reducible. In particular,
if $A>0$, then it is irreducible.

\begin{theorem}
Let $A=%
\llbracket
a_{j_{1}\cdots j_{k}}%
\rrbracket
\in\mathbb{R}^{n\times\dots\times n}$ be irreducible and $A\geq0$.
Then $A$ has a positive real $l^k$-eigenvalue with an
$l^k$-eigenvector $\mathbf{x}_{\ast}$ that may be chosen to have all
entries non-negative. In fact, $\mathbf{x}_{\ast}$ is unique and has
all entries positive.
\end{theorem}

\begin{proof}
Let $\mathbb{S}_{+}^{n}:=\{\mathbf{x}\in\mathbb{R}^{n}\mid\mathbf{x}%
\geq0,\lVert\mathbf{x}\rVert_{k}=1\}$. For any $\mathbf{x}\in\mathbb{S}%
_{+}^{n}$, we define%
\[
\mu(\mathbf{x}):=\inf\{\mu\in\mathbb{R}_{+}\mid A(I,\mathbf{x},\dots
,\mathbf{x})\leq\mu\mathbf{x}^{k-1}\}.
\]
Note that for $\mathbf{x}\geq0$, $\varphi_{k-1}(\mathbf{x})=\mathbf{x}^{k-1}$.
Since $\mathbb{S}_{+}^{n}$ is compact, there exists some $\mathbf{x}_{\ast}%
\in\mathbb{S}_{+}^{n}$ such that%
\[
\mu(\mathbf{x}_{\ast})=\inf\{\mu(\mathbf{x})\mid\mathbf{x}\in\mathbb{S}%
_{k}^{n}\}=:\mu_{\ast}.
\]
Clearly,%
\begin{equation}
A(I_{n},\mathbf{x}_{\ast},\dots,\mathbf{x}_{\ast})\leq\mu_{\ast}%
\mathbf{x}_{\ast}^{k-1}.\label{ineq}%
\end{equation}
We claim that $\mathbf{x}_{\ast}$ is a (mode-$1$) $l^{k}$-eigenvector of $A$,
ie.%
\[
A(I_{n},\mathbf{x}_{\ast},\dots,\mathbf{x}_{\ast})=\mu_{\ast}\mathbf{x}_{\ast
}^{k-1}.
\]
Suppose not. Then at least one of the relations in
\eqref{ineq}
must hold with strict inequality. However, not \textit{all} the relations in
\eqref{ineq}
can hold with strictly inequality since otherwise%
\begin{equation}
A(I_{n},\mathbf{x}_{\ast},\dots,\mathbf{x}_{\ast})<\mu_{\ast}\mathbf{x}_{\ast
}^{k-1}\label{sineq}%
\end{equation}
would contradict the definition of $\mu_{\ast}$ as an infimum. Without loss of
generality, we may assume that the first $m$ relations in
\eqref{ineq}
are the ones that hold with strict inequality and the remaining $n-m$
relations are the ones that hold with equality. We will write $\mathbf{x}%
_{\ast}=[\mathbf{x}_{0},\mathbf{x}_{1}]^{\intercal}$ with $\mathbf{x}_{0}%
\in\mathbb{R}^{m},\mathbf{x}_{1}\in\mathbb{R}^{n-m}$. By assumption, $A$
may be partitioned into blocks so that%
\begin{multline}
A_{00\cdots0}(I_{m},\mathbf{x}_{0},\dots,\mathbf{x}_{0},\mathbf{x}_{0})+\\
A_{00\cdots1}(I_{m},\mathbf{x}_{0},\dots,\mathbf{x}_{0},\mathbf{x}_{1}%
)+\dots+\\
A_{01\cdots1}(I_{m},\mathbf{x}_{1},\dots,\mathbf{x}_{1},\mathbf{x}_{1}%
)<\mu_{\ast}\mathbf{x}_{0}^{k-1},\label{top}%
\end{multline}%
\vspace*{-5ex}
\begin{multline}
A_{10\cdots0}(I_{n-m},\mathbf{x}_{0},\dots,\mathbf{x}_{0},\mathbf{x}_{0})+\\
A_{00\cdots1}(I_{n-m},\mathbf{x}_{0},\dots,\mathbf{x}_{0},\mathbf{x}%
_{1})+\dots+\\
A_{11\cdots1}(I_{n-m},\mathbf{x}_{1},\dots,\mathbf{x}_{1},\mathbf{x}_{1}%
)=\mu_{\ast}\mathbf{x}_{1}^{k-1}.\label{bot}%
\end{multline}
Note that $\mathbf{x}_{0}\neq\mathbf{0}$ since the \textsc{lhs} of
\eqref{top}
is non-negative. We will fix $\mathbf{x}_{1}$ and consider the following (vector-valued)
functions of $\mathbf{y}$:%
\begin{multline*}
F(\mathbf{y}) :=A_{00\cdots0}(I_{m},\mathbf{y},\dots,\mathbf{y}%
,\mathbf{y})+\\
A_{00\cdots1}(I_{m},\mathbf{y},\dots,\mathbf{y},\mathbf{x}%
_{1})+\dots+\\
A_{01\cdots1}(I_{m},\mathbf{x}_{1},\dots,\mathbf{x}_{1}%
,\mathbf{x}_{1})-\mu_{\ast}\mathbf{y}^{k-1},
\end{multline*}
\vspace*{-5ex}
\begin{multline*}
G(\mathbf{y}) :=A_{10\cdots0}(I_{n-m},\mathbf{y},\dots,\mathbf{y}%
,\mathbf{y})+\\
A_{00\cdots1}(I_{n-m},\mathbf{y},\dots,\mathbf{y},\mathbf{x}%
_{1})+\dots+\\
A_{11\cdots1}(I_{n-m},\mathbf{x}_{1},\dots,\mathbf{x}%
_{1},\mathbf{x}_{1})-\mu_{\ast}\mathbf{x}_{1}^{k-1}.
\end{multline*}
Let $f_{1},\dots,f_{m}$ and $g_{1},\dots,g_{n-m}$ be the component functions
of $F$ and $G$ respectively, ie.\ $f_{i}$'s and $g_{i}$'s are real-valued
functions of $\mathbf{y}\in\mathbb{R}^{m}$ such that
$F(\mathbf{y})=[f_{1}(\mathbf{y}),\dots,f_{m}(\mathbf{y})]^{\intercal}$ and
$G(\mathbf{y})=[g_{1}(\mathbf{y}),\dots,g_{n-m}(\mathbf{y})]^{\intercal}$.

By
\eqref{top}%
, we get $f_{i}(\mathbf{x}_{0})<0$ for $i=1,\dots,m$. Since $f_{i}$ is
continuous, there is a neighborhood $B(\mathbf{x}_{0},\delta
_{i})\subseteq\mathbb{R}^{m}$ such that
$f_{i}(\mathbf{y})<0$ for all $\mathbf{y}\in B(\mathbf{x}_{0},\delta_{i})$.
Let $\delta=\min\{\delta_{1},\dots,\delta_{m}\}$. Then
$F(\mathbf{y})<\mathbf{0}$ for all $\mathbf{y}\in B(\mathbf{x}%
_{0},\delta)$.

By
\eqref{bot}%
, we get $g_{j}(\mathbf{x}_{0})=0$ for $j=1,\dots,n-m$. Observe that if
$g_{j}$ is not identically $0$, then $g_{j}(\mathbf{y})=g_{j}(y_{1}%
,\dots,y_{m})$ is a non-constant multivariate polynomial function in the
variables $y_{1},\dots,y_{m}$. Furthermore, all coefficients of this
multivariate polynomial are non-negative since $A\geq0$. It is easy to see
that such a function must be `strictly monotone' in the following sense: if
$\mathbf{0}\leq\mathbf{y}\leq\mathbf{z}$ and $\mathbf{z}\neq\mathbf{y}$, then
$g_{j}(\mathbf{y})<g_{j}(\mathbf{z})$. So for
$\mathbf{0}\leq\mathbf{y}\leq\mathbf{x}_{0}$ and $\mathbf{y}\neq\mathbf{x}_{0}$,
we get $g_{j}
(\mathbf{y})<g_{j}(\mathbf{x}_{0})=0$. Since $A$ is irreducible,
$A_{10\cdots0}$ is non-zero and thus some $g_{j}$ is not identically $0$.

Let $\mathbf{y}_{0}$ be a point on the line joining $\mathbf{0}$ to
$\mathbf{x}_{0}$ within a distance $\delta$ of $\mathbf{x}_{0}$ and
$\mathbf{y}_{0} \ne \mathbf{x}_{0}$. Then with
$\mathbf{y}_{0}$ in place of $\mathbf{x}_{0}$, the $m$ strict inequalities in
\eqref{top}
are retained while at least one equality in
\eqref{bot}
will have become a strict inequality. Note that the homogeneity of
\eqref{top}
and
\eqref{bot}
allows us to scale $[\mathbf{y}_{0},\mathbf{x}_{1}]^{\intercal}$ to unit
$l^{k}$-norm without affecting the validity of the inequalities and
equalities. Thus we have obtained a solution with at least $m+1$ relations
in
\eqref{ineq}
being strict inequalities. Repeating the same arguments inductively, we can
eventually replace all the equalities in
\eqref{top}
with strict inequalities, leaving us with
\eqref{sineq}%
, a contradiction. [\emph{We defer the proof of uniqueness and
positivity of $\mathbf{x}_{\ast}$ to the full paper.}]
\end{proof}

A proposal to use the multilinear Perron-Frobenius theorem in the
ranking of linked objects may be found in \cite{L}. A symmetric
version of this result can be used to study hypergraphs \cite{DL}.

\bibliographystyle{IEEEbib}
\bibliography{camsap}

\end{document}